%added wright reference

\magnification=\magstep1
\topskip=60pt

\def\ie {{\it i.e.\ }}

\def\cf {\hbox{\it cf.\ }}

\def\Const {\hbox{Const}}

\def\Q {{\bf Q}}
\def\C {{\bf C}}
\def\Z {{\bf Z}}

\def\R {{\bf R}}
\def\vol{{\it vol}}

\def\sys{{\it sys}} 
\def\stsys{{\it stsys}} 
\def\mass{{\it mass}}

\def\ttf{{\tilde {\tilde f\phantom{i}}}}
\def\injrad{{\it InjRad}}
\def\SP{{\it SP}}
\def\NP{{\it NP}}
%\raggedright
\font\large=cmbx12 scaled \magstep1
\pageno=1
\normalbaselines 

\input epsf 		%%for the figure

%\headline{\hfill {\bf M. Katz} }
$\phantom{a}$
\vskip1cm
\medskip\noindent
\centerline{{\large Local calibration of mass and systolic geometry}}

\vskip1cm
\bigskip\noindent 

\centerline{{\large Mikhail
Katz}\footnote*{Supported by The Israel Science Foundation (grant no.\
620/00-10.0).  Partially supported by the Emmy Noether Research
Institute and the Minerva Foundation of Germany.  Paper visible at the
site http://www.cs.biu.ac.il/$\sim$katzmik/publications.html/ and will
appear in $\underline{\hbox{\it Geometric and Functional Analysis}}$
{\bf 12}, issue 3 (2002).}}

\vskip1cm

\bigskip\noindent {\bf Abstract:} We prove the simultaneous
$(k,n-k)$-systolic freedom, for a pair of adjacent integers $k<
{n\over 2}$, of a simply connected $n$-manifold $X$.  Our
construction, related to recent results of I. Babenko, is concentrated
in a neighborhood of suitable $k$-dimensional submanifolds of $X$.  We
employ calibration by differential forms supported in such
neighborhoods, to provide lower bounds for the $(n-k)$-systoles.
Meanwhile, the $k$-systoles are controlled from below by the
monotonicity formula combined with the bounded geometry of the
construction in a neighborhood of suitable $(n-k+1)$-dimensional
submanifolds, in spite of the vanishing of the global injectivity
radius.  The construction is geometric, with the algebraic topology
ingredient reduced to Poincare duality and Thom's theorem on
representing multiples of homology classes by submanifolds.  The
present result is different from the proof, in collaboration with
A. Suciu, and relying on rational homotopy theory, of the $k$-systolic
freedom of $X$.  Our results concerning systolic freedom contrast with
the existence of stable systolic inequalities, studied in joint work
with V.~Bangert.
\vskip1cm
\bigskip\noindent
{\large 1. Introduction}
\bigskip\noindent
Let $X$ be an orientable manifold of dimension $n$.  A choice of a
Riemannian metric $g$ on $X$ allows us to compute the total volume
$\vol_n(g)$, as well as the $k$-volumes of $k$-dimensional
submanifolds of $X$.  Given an integer homology class $\alpha\in H_k
(X,\Z)$, let $\vol_k(\alpha)=\inf_{x\in \alpha}\vol(x)$, where the
infimum is taken over Lipschitz cycles.  Here the volume of an integer
$k$-cycle $x=\sum_i r_i \sigma_i$ is $\vol(x)=\sum_i |r_i| \vol
(\sigma_i)$, while the volume of a Lipschitz singular chain $\sigma_i:
\Delta\to X$ is the integral over the standard $k$-simplex
$\Delta$, of the ``volume form'' of the pullback $\sigma_i^*(g)$.

We define the $k$-{\it systole} of $(X,g)$, denoted $\sys_k(g)$, as
the infimum of volumes of nonzero integer homology classes:
$$
\sys_k(g) =\inf_{\alpha\not= 0 }\vol_k(\alpha).
$$ 
Note that the $n$-systole coincides with the total volume:
$\sys_n(g)=\vol(g)$.  An alternative definition involves enlarging the
class of $k$-dimensional submanifolds to {\it rectifiable
$k$-currents}, so as to enable the solution of extremal problems, and
extending the $k$-volume to $k$-mass (\cf [Mo]).  Then the $k$-systole
is the least mass of a rectifiable $k$-current representing a
nonzero integer $k$-dimensional homology class.  For a leisurely
historical introduction to $k$-systoles, and for further references,
see M. Berger's survey article [Be2], or [KS1], sections 2 and 3.

M. Gromov asks the following question in his recent book [G3], p.~268.
Let ${\cal C}_l$ be the space of all left-invariant metrics on the
unitary group $U(d)$.  Let ${\cal S}_l: {\cal C}_l\rightarrow {\bf
R}^{d^2}$ be the `systolic' map, i.e.\ the map which associates to
every metric $g$, its vector of systolic invariants $sys_k(g)$,
$k=1,\ldots,\dim(U(d))=d^2$.  Gromov asked:

\medskip\noindent
{\bf Question 1.} What is the image of the map ${\cal S}_l$?
 
\medskip\noindent
More generally, let ${\cal C}_g(X)$ denote the space of arbitrary
Riemannian metrics $g$ on a given smooth $n$-manifold $X$.  Let ${\cal
S}_g: {\cal C}_g\rightarrow {\bf R}^n$ be the analogous `systolic'
map.

\medskip\noindent
{\bf Question 2.} What is the image of the map ${\cal S}_g$?

\medskip\noindent
Typically, the constraints on the image of ${\cal S}$ have been
expressed in terms of inequalities satisfied by the various systoles
or their ``stable'' analogues, \cf [G1], [H], [BanK], [K2].  The
absence of such inequalities, termed ``systolic freedom'', has been
exhibited in a number of recent works, starting with Gromov's example
of 1993.

Systolic freedom involving a single systole with regard to the total
volume has by now been well understood, at least in the case where the
relevant homology group is torsionfree.  Thus, even-dimensional
manifolds with middle-dimensional systole growing faster than the
volume were constructed in [BKS], generalizing [K1].  Manifolds whose
$k$-systole, for any given $2\leq k<n$, grows faster than the volume,
were constructed in [KS1], [KS2].

Systolic freedom of manifolds involving a pair of systoles of
complementary dimensions was studied in [BeK], [P], [BK].  Metrics $g$
on general $n$-dimensional polyhedra with the {\it product},
$\sys_k(g) \sys_{n-k}(g)$, growing faster than the total volume, were
constructed in [B2, B3], vastly generalizing the results of [BK] (\cf
section 2).

A shortcoming of such constructions is the fact that, as a by-product
of the construction, the injectivity radius tends to vanish, and
similarly for the systoles below dimension~$k$ 
%\vfill\eject\end
(\cf Remark 2.2 for
more details).  Thus, such constructions are limited to inequalities
involving only 2 systoles.

In the present work we attempt to tackle Question 2 above in a
situation where a larger number of systoles is involved, in the
following sense.

\medskip\noindent
{\bf Definition 1.0.}  Let $X$ be an $n$-dimensional smooth manifold,
and let $k$ be an integer satisfying $1\leq k < {n\over 2}$.  Then
$X$ is called {\it simultaneously $(k,n-k)$- and
$(k-1,n-k+1)$-systolically free} if
$$
\inf_g\vol(g)\left({1\over \sys_k(g)\sys_{n-k}(g)} +
{1\over \sys_{k-1}(g)\sys_{n-k+1}(g)}\right)=0,
\eqno{(1)}
$$
where the infimum is over all smooth metrics $g$ on $X$.

In other words, we show that the product of systoles of complementary
dimensions can be made large, compared to the volume, simultaneously
for {\it two} pairs of complementary dimensions, modulo mild
restrictions on the topology of the manifold.

We will adopt the convention that a systolic invariant defined over a
trivial homology group, is infinite.

\medskip\noindent
{\bf Theorem 1.1.}  Let $X$ be a simply connected closed
$n$-dimensional manifold.  Let $k<{n \over 2}$ be an integer.  Assume
that the group $H_{k-1}(X)$ is torsionfree.  Then $X$ admits
submanifolds $C_{k-1}$ and $C_k$ with $\dim (C_i)=i$, and a sequence
of metrics $(g_j)$ as $j\to \infty$, fixed outside a neighborhood of
the union $C_{k-1} \cup C_k\subset X$, which satisfy the following
three conditions:
\medskip\noindent
1.  We have a uniform lower bound for the $(k-1)$-systole and the
    $k$-systole: $\sys_{k-1} (g_j) \geq 1$ and $\sys_k (g_j) \geq 1$;
\medskip\noindent
2.  The total volume of $X$ is linear in $j$, namely $\vol(g_j)\leq
\mu(X) \sum_i b_i(X)\, j$, where $b_i$ are the Betti numbers, and
$\mu(X)$ is a constant depending only the topology of $X$ (\cf 5.1 and
5.3);
\medskip\noindent
3.  We have a quadratic lower bound for the complementary stable
    systoles: $\stsys_{n-k}(g_j)\geq j^2$ and $\stsys_{n-k+1}(g_j)\geq
    j^2$.

\medskip\noindent
Therefore, $X$ is simultaneously
$(k,n-k)$- and $(k-1,n-k+1)$-systolically free.

%\vfill\eject\end
\bigskip
Note that our method has the usual consequence in terms of systolic
freedom in a (single) pair of complementary dimensions, as in
Corollary 1.2 below.  This corollary follows from more general results
for arbitrary polyhedra proved in [B3], \cf section 2 below.

\medskip\noindent
{\bf Corollary 1.2.}  Let $X^n$ be orientable, let $k<{n\over 2}$, and
assume that the group $H_{n-k}(X)$ is torsion free.  Then $X$ is {\it
$(k,n-k)$-systolically free}.

\medskip\noindent
Indeed, the estimates above imply that the ratio
$${\vol_n(g_j) \over \sys_k(g_j) \sys_{n-k}(g_j)} \leq {\vol_n(g_j)
\over \sys_k(g_j) \stsys_{n-k}(g_j)} \leq {\mu(X) \sum_i b_i(X) j
\over 1\cdot j^2}= {\Const \over j}\to 0
\eqno{(2)}
$$ 
tends to zero.  Here we do not need the assumption of simple
connectivity of Theorem 1.1, which is only required for the
simultaneous construction, \cf section 6.

\medskip\noindent
Our theorem has the drawback of eliminating the interesting case of
dimension and codimension one.  The result below partly compensates
for it.

\medskip\noindent
{\bf Theorem 1.3.} (\cf 7.2 and 7.4) Let $X$ be an orientable
$n$-manifold, $n\geq 5$.  Then $X$ admits metrics which are
simultaneously $(1,n-1)$-free and $(2,n-2)$-free in the following two
cases: (a) the fundamental group $\pi_1(X)$ is free abelian; (b)
$H_1(X)=\Z$.

\medskip\noindent

We will place the present work in its context among other systolic
results in section~2.  The local construction is described in
section~3.  Bounded geometry is discussed in section~4.  Calibration
by differential forms in a neighborhood of a submanifold and the use
of the monotonicity formula for minimizing rectifiable currents appear
in section 5, which describes the main construction using
triangulations and pullback.  Theorems 1.1 is proved in section~6.  We
use the classification of surfaces and an analysis of systems of
disjoint loops on surfaces to prove Theorem 1.3 in section 7.

%\vfill\eject
\bigskip\noindent
{\large 2. Historical remarks and motivation}
\medskip\noindent
An early result on the $(k,n-k)$-freedom from [BK] (\cf Theorem 2.1
below) contained a modulo 4 restriction on $k$, as well as an
additional restriction of high connectivity.

I. Babenko [B2, B3] proved the systolic freedom in a pair of
complementary dimensions, $(k,n-k)$, of an arbitrary $n$-dimensional
polyhedron.  Babenko's proof [B3] involves polyhedral constructions
using products of spheres.  Babenko's result, which is valid in
particular for all manifolds regardless of their orientability, shows
that systolic freedom in a pair of complementary dimensions has
nothing to do with Poincare duality.

As compared to [B3], our theorem shows that such $(k,n-k)$-freedom can
be achieved simultaneously for a pair of distinct $k$, modulo the
stated assumptions on $X$.

More precisely, the following result follows from [BK], Lemma 5.1.

\medskip\noindent
{\bf Theorem 2.1} [BK].  Let $X$ be an orientable $n$-dimensional
manifold, and let $k < {n\over 2}$.  Assume that a multiple of every
$k$-dimensional homology class can be represented by a submanifold
$A\subset X$ with trivial normal bundle, which is either a sphere
$S^k$ or a product $A=B\times C$ where $C$ is a circle.  Then $X$
admits metrics $g_j$ with injectivity radius $\injrad(g_j)\geq 1$
uniformly bounded from below, sectional curvature uniformly bounded in
absolute value, while the total volume grows at most linearly:
$\vol(g_j) \leq j$, and the stable $(n-k)$-systole at least
quadratically: $\stsys_{n-k}(g_j)\geq \Const\, j^2$.  Thus,
$$
\injrad(g_j)^k \stsys_{n-k}(g_j) \geq j\, \Const \vol(g_j)
\eqno{(3)}
$$
where the constant, $\Const$, is independent of the metric.

\medskip\noindent
{\bf Remark 2.1.A: The monotonicity formula.}  The theorem implies, in
particular, the $(k, n-k)$-systolic freedom modulo suitable torsion
hypotheses.  This follows either from the coarea formula, or,
alternatively, from the monotonicity formula [Mo] applied to a
minimizer for $k$-dimensional homology, which gives a lower bound for
$\sys_k$ in terms of $(\injrad)^k$ (\cf equation $(13)$ below).

\medskip\noindent
A key ingredient in the proof of $(3)$ along the lines of [BK] in the
map
$$f: U_\epsilon A\to A\times D^{n-k},
\eqno{(4)}
$$ 
which identifies a neighborhood of a suitable representative, $A$, of
a homology class, with the product with the $(n-k)$-disk.

In contrast, the starting point of the argument in [B3] is the CW
complex $X/X^{(k-1)}$ obtained by collapsing the $(k-1)$-skeleton of
$X$ to a point.  Thus, one immediately loses control of the
injectivity radius and the systoles below dimension $k$.

The origin of the present paper was an attempt to understand the
argument of [B3].  It turned out that one can interpret this argument
as a local construction using a calibration, generalizing the map
$(4)$ from the construction of [BK] by introducing the degree 1 map
$f: U_\epsilon\to S^k\times D^{n-k}$ of formula $(7)$ below, where
$U_\epsilon=U_k$ is a tubular neighborhood of a submanifold
$C_k\subset X$ representing a $k$-dimensional homology class.

\medskip\noindent
{\bf Remark 2.2: Vanishing injectivity radius.}  The inclusion of
$C_k$ in $X$ may have nonzero image in the lower dimensional homology
groups.  Our construction relies on the pullback of metrics by the map
$f$ of $(7)$, which annihilates the lower dimensional homology of
$C_k$.  In particular, the resulting metrics may have a vanishing
$(k-1)$-systole.  The author is grateful to I. Babenko for pointing
out this drop in the lower systoles, analogous to the collapse of the
$(k-1)$-skeleton in the construction of [B3].

Note that $(k,n-k)$-systolically free metrics on $X$ are obtained in
[B3] by pullback by a map $f: X \to W$, where $W$ a polyhedron
obtained from $X/X^{(k-1)}$ by attaching copies of the product of
spheres, $S^k\times S^{n-k}$.  Here [B3] uses $(k,n-k)$-free metrics
from [BK] on the product of spheres.  Meanwhile, the construction of
the $(n-k)$-free metrics on the product of spheres of [BK] is local in
a neighborhood of the factor $S^k$.  Thus, ultimately, the
construction of [B3] relies on the construction of $(k,n-k)$-freedom
in a neighborhood of a specific $k$-cycle in $X$.  In the case when
$X$ is a manifold, such a specific cycle is the inverse image of a
suitable copy of $S^k$, or, more precisely, the inverse image of a
regular value of the composition of $f$ with the projection to the
second factor $S^{n-k}$.  It is the latter observation that originally
motivated the author to undertake the present work.

Placing the local nature of the construction in the forefront allows
us to prove {\it simultaneous} freedom in two adjacent pairs of
complementary dimensions, while remaining in the category of
manifolds.  The input from algebraic topology is reduced to Poincare
duality and Thom's theorem, as follows.

To protect the $(k-1)$-systole from collapse, we sharpen the
construction of [BK] to include control over $(n-k+1)$-dimensional
submanifolds $B_{n-k+1}\subset X$.  More precisely, the geometry in a
neighborhood of $B_{n-k+1}$ remains uniformly bounded.  This provides
a uniform lower bound for the volume of a minimizer for an infinite
order integer $(k-1)$-dimensional homology class, by Poincare duality
and the monotonicity formula applied in a neighborhood of $B_{n-k+1}$.

\medskip\noindent
{\bf Remark 2.3: Relation to $k$-systolic freedom.}  If we rescale the
metric $g_j$ of Theorem~1.1 to {\it unit} total volume, then the
$(n-k)$-systoles above the middle dimension become arbitrarily large,
while the $k$-systoles below the middle dimension tend to 0.  In other
words, these metrics are ``$(n-k)$-systolically free'', but not
``$k$-systolically free''.

Can the systoles below the middle dimension also be made arbitrarily
large, compared to the total volume of $X$?  The answer is typically
negative for $k=1$.  Namely, M. Gromov proved the following result in
[G1].  Let $\sys \pi_1(g)$ be the homotopy 1-systole of $(X,g)$, \ie
the length of the shortest noncontractible loop.  Assume $X$ is
aspherical, or, more generally, ``essential'', \cf [G1].  Then we have
the inequality $$ \sys \pi_1(g)\leq C_n \vol_n(g)^{{1\over n}}, $$
where the constant $C_n$ depends only on the dimension.  

On the other hand, the answer is affirmative for $k\geq 2$, for an
arbitrary manifold $X$ with torsion free $k$-dimensional homology, as
shown by the author in collaboration with A.~Suciu [KS1, KS2].
Compared to the $(k,n-k)$-freedom, such results seem to be harder to
obtain in a purely geometric way.  Thus, they seem to require
``classifying space'' type arguments from algebraic topology.  The
proof of [KS1, KS2] uses rational homotopy theory.  In contrast, the
present paper uses direct geometric constructions.

\medskip\noindent
{\bf Question 2.4.} Can $k$-systolic freedom be attained
simultaneously for more than a single dimension $k\leq {n \over 2}$?

\medskip\noindent
Our theorem can be usefully compared to the following result of
J. Hebda [H], \cf [BanK].

Let $H_k(X,\Z)_\R\subset H_k(X,\R)$ be the maximal lattice obtained as
the image of integer homology.  Denote by $\alpha_\R\in H_k(X,\Z)_\R$
the image of the class $\alpha\in H_k(X,\Z)$.  We set $$
\mass_k(\alpha_\R)=\lim_{i\to \infty} {1\over i} \vol(i\alpha).
%\eqno{(?)}
$$ We define the stable $k$-systole, denoted $\stsys_k(g)$, by
minimizing mass over nonzero elements in the integer lattice: $$
\stsys_k(g)=\inf_{\alpha_\R\not= 0 \in H_k(X,\Z)_\R}\mass_k(\alpha_\R).
$$
Alternatively, the stable $k$-systole is the least mass of a
minimizing {\it normal} $k$-current representing a nonzero element of
the lattice $H_k(X,\Z)_\R$.

\medskip\noindent
{\bf Theorem 2.5.} [H] Let $X$ be a compact orientable manifold of
dimension $n$.  Let $k< n$ and assume that the $k$-th Betti number is
positive, $b_k(X)>0$.  Then every metric $g$ on $X$ satisfies the
inequality
$$\stsys_k(g) \stsys_{n-k}(g)\leq C(n)b_k(X) \vol(g),
$$ 
where $C(n)$ depends only on the dimension $n$ of $X$.

\medskip\noindent
Gromov's paper [G1, section 7.4] contains general results of this
type; see also [BanK].  The combination of our Theorem 1.1 with
Hebda's Theorem 2.5 shows that the ratio
$${\sys_k(g_j)\over \stsys_k(g_j)} \geq O(j)
\eqno{(5)}
$$ tends to infinity.  Thus our theorem can be viewed as a systematic
way of constructing minimizing rectifiable currents which are far from
minimizing as normal currents.

The subject of systolic freedom has received renewed attention
recently in connection with the theoretical work in the context of
quantum computers.  Thus, M.~Freedman [Fr] proved that the 3-manifold
$S^1\times S^2$ admits (1,2)-systolically free metrics {\it even when
we allow nonorientable surfaces} to compete in the definition of the
2-systole.  Note that the case of orientable surfaces is easier and
was proved in [BeK] and generalized in [P] and [BK].  The general
study of systoles was pioneered by M. Berger [Be1]; see Gromov's 1999
book [G3] for an overview.

%\vfill\eject
\bigskip\noindent
{\large 3. The local ``two circle'' construction}
\medskip\noindent
Our construction is local in a neighborhood $U_k$ of a $k$-dimensional
submanifold $C_k\subset X$, specified in the lemma below.  We will
refer to it as the ``two circle'' construction, because it involves
the splitting off of a pair of circles, denoted $C$ and $S^1$.  The
circles are split off, respectively, of the $k$-dimensional class and
of the $(n-k)$-dimensional class, as in Proposition 3.3 below.

Let $(c_i)$ be an integer basis for a maximal lattice in $H_k(X,\Z)$.
The following lemma is well known [T] (\cf Theorem 6.2 below).
\medskip\noindent
{\bf Lemma 3.1.}  Let $k<{n\over 2}$.  There exists a submanifold
$C_k\subset X$, and an integer $\lambda\in \Z$, satisfying the
following two conditions:
\medskip\noindent
(a) the connected components $C_{k,i}$ of $C_k$ represent a fixed
integer multiple of the basis $c_i$, namely $[C_{k,i}] =
\lambda c_i$ for all $i$;
\medskip\noindent
(b) the normal bundle of $C_k$ in $X$ is trivial.
\medskip\noindent
In view of condition (b), a neighborhood $U_k$ of $C_k$ is
diffeomorphic to a product $C_k\times D^{n-k}$ with a ball $D^{n-k}$.
Let $(m_i)$ be a basis for $H_{n-k}(X)$ modulo torsion, dual to the
basis $(c_i)$.  Representative fibers $D^{n-k}$ from the connected
components of the fibration $U_k\to C_k$ may be ``closed up'', by
Poincare duality and hypothesis (a), to $(n-k)$-cycles $M_{n-k,i}
\subset X$ representing the elements of the basis $(m_i)$.  Thus the
algebraic intersection numbers with the components of the $k$-cycle
$C_k$ satisfy the relation
$$C_{k,i} \cdot M_{n-k, \ell}=\lambda \delta_{i\ell}.
\eqno{(6)}
$$ 
We have the excision isomorphism $r: H_{n-k}(U_{k,i}, \partial
U_{k,i}) \to H_{n-k}(X, X\setminus U)$ relating the classes $[M_i] =
r([D^{n-k}])$, where
$$[D^{n-k}]\in H_{n-k}(U_{k,i}, \partial U_{k,i})
$$
is a relative homology class.  

\medskip\noindent
{\bf Definition 3.2.}  Consider the map
$$f_k: U_k\to S^k\times D^{n-k},
\eqno{(7)}
$$ where $S^k$ is the $k$-sphere, defined by sending each connected
component of $C_k$ to $S^k$ by a degree 1 map (or any nonzero degree),
where $f_k$ is the identity on the second factor $D^{n-k}$.  

\medskip\noindent
If $\alpha\subset U_k$ is a $k$-cycle representing a nonzero homology
class, then the class of $f(\alpha)\in H_k(S^k\times D^{n-k})$ is also
nonzero.  Indeed, consider the generator $[C_k]\in H_k(U_k)=\Z$, and
let $[\alpha] = a[C_k]$, where $a\not= 0$.  Then the image under 
$f=f_k$ is
$$
[f(\alpha)] = f_*([\alpha])= f_*(a[C_k])=a f_*([C_k])= a \deg(f)
[S^k]
\not= 0
$$
in the group $H_k(S^k \times D^{n-k})=\Z$.

\medskip\noindent
Our target space, $S^k\times D^{n-k}$, of the map $f$ plays a key role
in the construction.  Namely, we will pull back ``systole-rich''
metrics from the target to the source, using a simplicial
approximation of $f$ with respect to suitable triangulations of source
and target.  For the purpose of proving ``simultaneous'' freedom, we
will need a simplicial approximation which remains a diffeomorphism in
certain neighborhoods where $f$ itself is a diffeomorphism.

\medskip\noindent
{\bf Proposition 3.3 (The two-circle construction).} Let $C\subset
S^k$ be a distinguished circle.  Then the manifold $Y= S^k\times
D^{n-k}$ admits a submanifold $S^{k-1}\times C$, where the
distinguished circle occurs as a copy of the second factor, with the
following property.  A neighborhood $Y'=S^{k-1}\times C\times D^{n-k}$
of $S^{k-1}\times C$ contains a hypersurface
$$\Sigma= S^{k-1}\times C\times S^1 \times K,
$$ 
with a neighborhood $Y'' = \Sigma \times I$, where $S^1$ is a circle,
while $\dim(K)=n-k-2$, and furthermore:

\medskip\noindent
(i) we have $h([S^{k-1}\times C])= [S^k]$ for the inclusion
homomorphism $h: H_k(\Sigma)\to H_k(Y)$;

\medskip\noindent
(ii) we have $h'([S^1\times K \times I])= [D^{n-k}]$ for the excision
isomorphism
$$h': H_{n-k}(Y'',\partial Y'')\to H_{n-k}(Y,Y\setminus int(Y'')),
$$
where $int$ denotes interior.

\medskip\noindent
{\it Proof.}  Let $C\subset S^k$ be a circle, and consider the class
$[S^k]\in H_k(Y)$.  It can clearly be represented by an imbedded
submanifold $S^{k-1}\times C\subset Y$ with trivial normal bundle.  A
tubular neighborhood $Y'$ of $S^{k-1}\times C$ in $Y$ is diffeomorphic
to the product
$$
Y'=S^{k-1}\times C\times D^{n-k}.
\eqno{(8)}
$$  
All the $j$-dependent constructions will take place inside $Y'$.  Now
consider a codimension 2 submanifold $K\subset D^{n-k}$ (for example,
an $(n-k-2)$-sphere), still with trivial normal bundle.  A
neighborhood of $K$ in $D^{n-k}$ is diffeomorphic to $D^2\times K$,
with boundary $S^1\times K$.  Thus a neighborhood of $S^{k-1}\times
C\times K$ in $Y'$ is diffeomorphic to $S^{k-1}\times C\times D^2
\times K$, with boundary denoted 
$$\Sigma= S^{k-1}\times C\times S^1 \times K.$$ Now a tubular
neighborhood $Y''$ of $\Sigma$ in $Y'$ is diffeomorphic to 
the product
$$
Y''=\Sigma\times I,
$$ 
where $I$ is an interval.  All the $j$-dependent
constructions will take place inside $Y''=\Sigma\times I$.

\bigskip\noindent
{\large 4. Bounding geometry near a submanifold while $\injrad=0$}

\medskip\noindent
We continue with the notation of the previous section.  Our goal is to
construct systolically free metrics while retaining a uniform bound on
the geometry in a neighborhood of a suitable submanfold $Z\subset
Y=S^k\times D^{n-k}$.

\medskip\noindent
{\bf Definition 4.1.}  Denote by $Z\subset Y''\subset Y$ the
submanifold $Z=C\times S^1\times K\times I$ with boundary, so that
$Y''=S^{k-1}\times Z$.

\medskip\noindent
{\bf Definition 4.2.}  We will say that a family $g_j$ of metrics on
a manifold with boundary is of ``bounded geometry'' if it has the
following three properties:

(a) the sectional curvature $K$ is bounded, $|K|\leq 1$, uniformly in
$j$; 

(b) the metric is constant (independent of $j$) in a unit neighborhood
of the boundary;

(c) the injectivity radius $\iota$ satisfies $\iota_x(g_j)\geq 1$, for
all points $x$ at least a unit distance away from the boundary.

\medskip\noindent
{\bf Definition 4.4.}  We will say that a $p$-form $w$ is {\it
calibrating} if $w(X_1,\ldots,X_p)\leq 1$ for all $p$-tuples of unit
vectors $X_i$ (\cf [G2], 4.$A_1$).

\medskip\noindent
By the natural pairing $\langle w,\phantom{a} \rangle$ between
homology and cohomology, an $(n-k)$-form $w$ with compact support can
be integrated over $D^{n-k}\subset Y$, viewed as a relative cycle
defining a class in the group $H_{n-k}(Y,Y\setminus Y'')=
H_{n-k}(Y'',\partial Y'') $.

\medskip\noindent
{\bf Proposition 4.5.} The manifold $Y''\subset Y= S^k\times D^{n-k}$,
where
$$Y''=\Sigma\times I = S^{k-1}\times C\times S^1 \times K \times I
$$
admits metrics $\tilde g_j$, together with closed $(n-k)$-forms $w_j$
compactly supported in the interior of $Y''$, with the following four
properties:
\medskip\noindent
(i) the metrics have ``bounded geometry'' in the sense of 4.2;

\medskip\noindent
(ii) the total volume is linear in $j$, \ie $\vol(\tilde g_j)\leq j$;
\medskip\noindent
(iii) we have the ``area'' lower bound $\langle w_j, D^{n-k} \rangle
\geq j^2$;
\medskip\noindent
(iv) the $(n-k+1)$-volume of the submanifold $Z= C\times S^1\times
K\times I\subset Y''$ is linear in $j$.
\medskip\noindent
{\it Proof.} 
The starting point of the construction is the fundamental
domain for the manifold 
$$H/G(j),
$$ 
where $H\subset SL(3,\R)$ is the Heisenberg group of unipotent
matrices, while $G(j)\subset SL(3,\Z)$ is the subgroup consisting of
matrices
$$
\left[ 
\matrix{1&x&z\cr 0&1&y\cr
0&0&1\cr}\right];\ x,y,z\in\Z
$$
such that $x$ is congruent to $0$ modulo $j$ (\cf [G3, p.~88] and [G2,
section 3.$C_6$]).

Suitable metrics $g_j$, and calibrating 2-forms $\psi_j$, on the
3-manifold with boundary $T^2\times I=C\times S^1\times I$ were
constructed in [BK] (see also [KS1, Appendix A]).  These are related
to F.~Almgren's example on $S^1\times S^2$ (\cf [Fe2], p.~397).

Here one may equip the manifold $Y= T^2\times I$ with metrics $g_j$
defined by
$$
g_j(x,y,z)= h(\hat x) (y,z) + dx^2,
$$
where $x\in I=[0,2j]$, $\hat x= \min (x,2j-x)$, while $T^2$ is the
quotient of the $(y,z)$-plane by the integer lattice, and  the formula
$$h(x) (y,z)
= (dz-xdy)^2 + dy^2
$$ 
defines a metric on the $2$-torus $T^2\times
\{x\}$.  Here the $z$-axis parametrizes the circle $C$ which has unit
length, while the $y$-axis parametrizes the circle $S^1$ of length
$\sqrt{1+(\hat x)^2}$.  Note that the two circles, $C$ and $S^1$, play
very different roles in this seminal ``two-circle'' construction.

The calibrating form is $\psi_j=(1+x^2)^{-{1\over 2}}*dz$, where $*$
is the Hodge star operator of the metric $g_j$.  For details, see [BK]
and [KS1], Appendix A.

The 1-systole $\sys_1(g_j)$ is uniformly bounded from below, while the
stable 1-systole tends to zero at the rate of ${1\over j}$.  The
volume grows linearly in $j$.

Meanwhile, the area of the surface $S^1\times I\subset T^2\times I$
grows quadratically in $j$ by calibration by $\psi_j$.  Such growth is
reflected in the fact that the minimizing surface $S^1\times I$
accumulates on itself, in the sense that distinct sheets of the
surface squeeze together to within distance on the order of ${1\over
j}$, comparable to the stable 1-systole.  In contrast, the injectivity
radius remains uniformly bounded from below.

Here a basic building block is a fundamental domain for the standard
compact nilmanifold of the Heisenberg group with its left invariant
metric.

We imbed $(T^2\times I, g_j)$ in $T^2\times K\times I$, by
incorporating a fixed (independent of $j$) metric on $K$ as a direct
summand, and let $w_j= \psi_j \wedge d_{vol_K}$ be the $(n-k)$-form
obtained by exterior product with the volume form of $K$.  Finally, we
pull the form $w_j$ back to $Y''=\Sigma\times I$ by the coordinate
projection $S^{k-1}\times T^2\times K\times I \to T^2\times K \times
I$, to obtain the desired form on $Y''$.  Here the metric of the
factor $S^{k-1}$ is a direct summand of the metric on
$Y''=S^{k-1}\times (T^2\times K\times I)$.

The quadratic lower bound for an $(n-k)$-cycle $\alpha$ follows, as in
[BK] and [KS1], by integrating the calibrating form $w_j$ over $\alpha$
in the connected component of $U_k$ containing a $C_{k,i}$ which has
nonzero algebraic intersection with $\alpha$.

Finally, the linear upper bound for the total volume as well as the
$(n-k+1)$-volume of the manifold $Z=C\times S^1\times K\times I$
follows from the linear upper bound for the volume of the metrics on
the 3-manifold $T^2\times I$, together with the fact that the metrics
$\tilde g_j$ on $Y''$ are a direct sum with a fixed metric on the
factors $K$ and $S^{k-1}$.

\medskip\noindent
{\bf Lemma 4.6.}  The $k$-systole of the metrics $\tilde g_j$ on $Y''$
is uniformly bounded from below.

\medskip\noindent
{\it Proof.} Let $\alpha$ be a $k$-cycle in $Y''$ representing a
nonzero multiple of the class $[S^{k-1}\times C]$.  Recall that
$Y''=S^{k-1}\times Z$ where $Z=C\times S^1\times K\times I$.  Consider
the projection to the first factor $p:Y''\to S^{k-1}$.  Then $p$ is a
Riemannian submersion by construction.  By the coarea (Eilenberg's)
formula [Mo], p.~31, we have
$$\vol_k(\alpha)\geq \int_{S^{k-1}} \alpha\cap p^{-1}(x) dx,
\eqno{(9)}
$$
where the intersection $\alpha\cap p^{-1}(x)$ with a typical fiber $Z$
is a loop in $Z$ which is not nullhomologous.  Hence $\sys_k(Y'')\geq
\vol_k(\alpha) \geq sys_1(Z) \vol_{k-1}(S^{k-1})$, providing the desired
lower bound in view of the bounded geometry of $Y''$ and $Z$.  Note
that the {\it stable} 1-systole of $Z$ tends to zero as $j$ increases.

However, we will provide an alternative argument, as well, which will
work even in a situation where the projection $p$ to $S^{k-1}$ is not
available, and the geometry far away from a specific copy of $Z$ may
not be bounded.  This is due to the collapse which occurs under the
map $f_k$ of formula $(7)$ above.  Therefore we would like to argue in
a fixed neighborhood of the submanifold $Z=C\times S^1 \times K\times
I$.  Indeed, the bounded geometry in a neighborhood of $Z$ allows us
to apply the monotonicity formula, centered at a point of $Z$, to a
minimizer in the class $[\alpha]$, immediately yielding the desired
lower bound.

%\vfill\eject
\bigskip\noindent
{\large 5. Triangulations, pullback metrics, and calibrations}
\medskip\noindent
In this section we will prove the $(k,n-k)$-systolic freedom for a
single $k$, preparing the ground for simultaneity in the next section.

The idea is to construct first the metrics $\tilde g_j$ on the product
$S^k\times D^{n-k}$, pull them back to $X$ by the map $f_k$ of formula
$(7)$, and then apply calibration by the pullbacks of the form $w_j$
of Proposition 4.5 to obtain the quadratic lower bound for the
$(n-k)$-systole.

By the work of I. Babenko [B1], the map $f_k: U_k\to S^k\times
D^{n-k}$ of formula $(7)$ may be replaced by a map $\tilde f$ which
has the following property with respect to suitable triangulations of
its domain and target: on each simplex, $\tilde f$ is either a
diffeomorphism onto its image, or the collapse onto a wall of positive
codimension.  Moreover, we will use the following theorem of
A. H. Wright (already exploited in [B1], Theorem 8.1).  Recall that a
map is called {\it monotone} if the inverse image of every point in
the range is a compact connected subset of the domain.

\medskip\noindent
{\bf Theorem 5.1.} ([W], Theorem 7.3.)  Let $M\sp n$ and $N\sp n$ be
closed piecewise linear manifolds. Let $f\colon M\sp n\rightarrow N\sp
n$ be a continuous map of absolute degree one. Then $f$ is homotopic
to a piecewise linear monotone map. 

\medskip
Note that if $M\sp n$ and $N\sp n$ are orientable, and if $f$ has
degree $±1$, then $f$ has absolute degree one.  In the terminology of
section 3, the map from the connected component $C_{k,i}$ to $S^k$ may
be chosen to be monotone by Wright's theorem.  Hence there is exactly
one $k$-dimensional simplex mapping diffeomorphically to a
$k$-dimensional simplex in the target $S^k$.  By Cartesian product
with the disk, we may assume that the map $\tilde f$, defined on each
connected component of the neighborhood of the submanifold $C_k$ of
Lemma 3.1, has the property that each top dimensional simplex in the
range has a unique simplex as its inverse image.

\medskip\noindent
The map $\tilde f$ can be slightly perturbed to $\ttf$ to make it
smooth, so the pullbacks by $\ttf$ are well-defined, while not
significantly affecting the volumes.  Consider the pulled-back
positive (symmetric) 2-forms $\ttf^*(\tilde g_j)$.  By construction,
its total volume is nearly equal to that of $\tilde g_j$.  The form
$\ttf^*(\tilde g_j)$ may not be definite, but a small multiple of a
fixed metric on $X$ can always be added on to make it positive
definite, in the end of the construction, without significantly
increasing the total volume, and certainly not decreasing the systolic
invariants.

\medskip\noindent
{\bf Proposition 5.2.}  Let $X$ be an orientable $n$-manifold (not
necessarily simply connected).  Let $k < {n \over 2}$.  Then the
metrics $g_j$ on the neighborhoods $U_k$ of Lemma 3.1, constructed by
pullback, can be extended from $U_k$ to $X$ so as to satisfy a uniform
lower bound for the $k$-systole, a quadratic (in $j$) lower bound for
the stable $(n-k)$-systole, and a linear (in $j$) upper bound for the
total volume.

\medskip\noindent
{\it Proof.}  Outside the neighborhood $U_k$, the metric $g=g_j$ on
$X$ is chosen fixed.  To patch together the metric $\ttf^*(\tilde
g_j)$ on $U_k$ and the metric $g$ on $X\setminus U_k$, we consider the
boundary $\partial U_k = C_k \times S^{n-k-1} \subset X$.  Its tubular
neighborhood in $X$ is of the form
$$C_k \times S^{n-k-1} \times I \subset  X.
$$ 
We use a partition of unity along the ``cylinder'' $C_k \times
S^{n-k-1} \times I$ (all choices independent of $j$).  Let $\alpha$ be
an integer $k$-cycle, which represents a nontrivial homology class.
We replace $\alpha$ by the minimizing rectifiable current in its
homology class.

\medskip\noindent
{\bf Remark 5.3.}  If $\alpha$ lies in $U_k$, then a uniform lower
bound
$$
\vol_k(\alpha)\geq \mu(X)
$$
for $\vol_k(\alpha)$ follows from the second argument in the proof of
Lemma 4.6.  If $\alpha$ ventures outside $U_k$, then such a lower
bound for its volume follows from the monotonicity formula centered at
a point of $\alpha\cap (X\setminus U_k)$, together with the fact that
the metric on $X\setminus U_k$ is fixed.

\medskip\noindent
{\bf Remark 5.4.}  Note that there is no lower bound for the
injectivity radius of $\ttf^* (\tilde g_j)$ , since the map $\ttf$
may collapse top-dimensional cells.

\medskip\noindent
Next, the total volume of $(X, g_j)$ is dominated by $j$ times the
number of connected components of $U_k$ by Theorem 5.1 and property
(ii) of Proposition 4.5.  This number is controlled by the Betti
number $b_k$ in view of Lemma 3.1(a), proving the upper bound for the
volume.

Now the pullback form $\ttf^*(w_j)$ on $U_k$ may be extended by zero
to give a calibrating form on all of $X$.  Let $M$ be a Lipschitz
$(n-k)$-cycle representing an infinite order class $[M]\in
H_{n-k}(X)$.  We decompose it with respect to the basis $(m_i)$:
$$[M]= T+ \sum_i \alpha_i m_i,
$$ 
where $T \in H_{n-k}(X)$ is torsion and $\alpha_i\in\Z$.  Choose a
nonzero index $\alpha_{i_0}\not=0$.  Let $w_{j,i_0}$ be the pullback
of $w_j$ supported in the connected component $U_{k,i_0}$.  Thus it is
the extension by zero of the pullback of $w_j$ by the restriction of
the map $f$ to this component.  Then by Stokes' theorem,
$$\int_M w_{j,i_0} = \int_{M\cap U_{n-k,i_0}} w_{j,i_0} = \alpha_{i_0}
\int_{M_{n-k,i_0}} \ttf^*(w_j) = \alpha_{i_0}\int_{D^{n-k}} w_j,
$$ where the $M_{n-k,i}$ are the cycles from formula (6).  Since $w_j$
is a calibrating form and $\ttf$ is an isometry, we have
$$\vol_{n-k}(M)\geq
|\alpha_{i_0}|\int_{D^{n-k}} w_j \geq j^2
\eqno{(10)}
$$ by property (iii) of Proposition 4.5.  Since the lower bound is
obtained by integration of differential forms, it holds for an
arbitrary normal current, proving the proposition.
\bigskip\noindent
{\large 6. Proof of simultaneous freedom}

\medskip\noindent
{\bf Theorem 6.1.} Let $X$ be a simply connected closed manifold of
dimension $n$.  Let $k< {n \over 2}$ and assume that $H_{k-1}(X)$ is
torsion free.  Then $X$ is simultaneously $(k,n-k)$-free and $(k-1,
n-k+1)$-free, in the precise sense stated in Theorem 1.1.

\medskip\noindent
{\it Proof.}  The proof breaks up into 5 steps.  Simple connectivity
is used in step 3, while torsion freeness is used in step 4.

\medskip\noindent
{\bf Step 1: Background metric and injectivity radius.}  Consider a
fixed smooth background Riemannian metric $g(0)=g_j(0)$ on $X$, with
positive fixed injectivity radius.

\medskip\noindent
{\bf Step 2: The $(k,n-k)$-freedom.}  Choose a rational basis for
$H_k(X)$, representable by manifolds $C_{k,i}$ with trivial normal
bundles.  Modify the background metric in a neighborhood $U_k$ of each
$C_{k,i}$, to construct $(k,n-k)$-free metrics $g_j(1)$ on $X$ as in
Proposition 5.2.  Recall that we have the submanifold $Z= C\times
S^1\times K\times I$ of $Y''=S^{k-1} \times Z$ which occurred in the
construction of the metrics.  The metrics $g_j(1)$ have the following
four properties:

\medskip\noindent
(a) the
$k$-systole is uniformly bounded from below: $\sys_k(g_j(1)) > 1$; 

\medskip\noindent
(b) the $(n-k)$-systole grows quadratically: $\sys_{n-k}(g_j(1)) =
j^2$;

\medskip\noindent
(c) the volume grows at most linearly: $\vol_n(g_j(1)) < j$.

\medskip\noindent
(d) the volume of the submanifold $Z$ grows linearly in $j$.

\medskip\noindent
Note that the construction of the $(k,n-k)$-free metrics results in
metrics $g_j(1)$ with zero injectivity radius.  Notwithstanding, the
uniform lower bound for the $k$-systole is ensured by Proposition 5.2,
in view of the bounded geometry in a neighborhood of $Z$.

\medskip\noindent
{\bf Step 3: Choice of the $(n-k+1)$-dimensional submanifold $B$.}
Unlike the uniform lower bound for the $k$-systole, in the wake of
step 2, we have no control over either the $(k-1)$-systole or the
injectivity radius, which may both be zero if we work with positive,
possibly nondefinite, quadratic forms $g_j(1)$.  To repair this and
ensure a lower bound for the $(k-1)$-systole, we need to modify the
construction as follows.

Choose a rational basis for $H_{n-k+1}(X)$ which can be represented by
$b_{k-1}$ imbedded submanifolds $B= B_{(n-k+1),i}$.  See
Sh.~Weinberger's comments below in Theorem 6.2, concerning the
existence of the submanifolds $B$.

We will modify the construction so as to have precise control over the
metric in a neighborhood of $B$.  By Poincare duality, every
$(k-1)$-cycle $\alpha$ representing a class of infinite order must
meet one of the submanifolds $B$.  This will yield a lower bound for
the $(k-1)$-volume of $\alpha$.

If $B$ lies outside the neighborhood $U_k$ of $C_k$, where the
$j$-dependent construction took place in step 2, then, by step 1,
there is a lower bound for the injectivity radius in a fixed
neighborhood of $B$, yielding the desired lower bound for a
$(k-1)$-cycle with nonzero algebraic intersection with $B$.

Now assume that the connected imbedded submanifold $B_{n-k+1,i}$ meets
$C_k$.  We place them in transverse position.  Then 
$$
\dim(B_{n-k+1,i} \cap C_k)= \dim(B_{n-k+1,i})+\dim(C_k)-n =1,
$$ 
hence the intersection $B_{n-k+1,i}\cap C_k$ is a disjoint union of
circles.

If $k=2$, then the $(k-1,n-k+1)=(1,n-1)$-freedom is not an issue,
since $X$ is assumed simply connected (\cf section 7 for a result on
non-simply connected manifolds in the case $k=2$).  Therefore we may
assume that $k\geq 3$.  Then by transversality we may assume that the
triple intersections
$$
B_{n-k+1,i}\cap B_{n-k+1,j}\cap C_k= \emptyset
$$ 
are empty (\cf formula $(14)$ of section 7).  Thus the full
intersection $(\cup_i B_{n-k+1,i})\cap C_k$ is also a disjoint
collection of circles, denoted $C_\alpha$:
$$
(\cup_i B_{n-k+1,i})\cap C_k=\cup_\alpha C_\alpha.
$$
A tubular neighborhood of $C_\alpha\subset B$ is of the form
$$
C_\alpha \times D^{n-k}\subset B.
\eqno{(11)}
$$
Recall that $X$ is assumed simply connected.  Therefore, after a
suitable surgery, we can assume that the manifold $C_k$ is simply
connected, as well.  Hence, the circle $C_\alpha$ is contractible in
$C_k$.  Choose a neighborhood $A_\alpha$ of an imbedded disk which
bounds $C_\alpha \subset C_k$, such that $A_\alpha$ is diffeomorphic
to a $k$-ball.  We choose such neighborhoods disjoint for distinct
circles $C_\alpha$.

Now choose a fixed circle $C\subset S^k$.  Let $f: C_k\to S^k$ be a
map with the following three properties.

1) $f$ maps the interior of each $A_\alpha$ diffeomorphically onto the
complement of a fixed basepoint in $S^k$;

2) $f$ collapses the complement $C_k\setminus (\cup_\alpha A_\alpha)$
   to the basepoint;

3) Each circle $C_\alpha$ is mapped diffeomorphically to the fixed
   circle $C\subset S^k$.

\medskip\noindent
Note that the degree of $f$ is the total number of circles $C_\alpha$.
We may use our map $f$ to define the map $f_k: U_k \to Y= S^k\times
D^{n-k}$, in place of the degree one map of formula~$(7)$.  The
diffeomorphism property 1) above allows us to argue as if a
neighborhood of a circle $C_\alpha\subset B$ is actually contained in
the standard manifold $Y$, and apply transversality techniques there
to achieve the identification of formula $(12)$ below.

We can assume that the image $f(C_\alpha)$ of the circle $C_\alpha
\subset C_k$ coincides with a  circle $C$ which is a copy
of the second factor in the submanifold $S^{k-1} \times C \subset Y$
homologous to $S^k$, and that the product structure $(11)$ in $B$
coincides with the last two components of decomposition $(8)$ of
$Y'=S^{k-1}\times C\times D^{n-k}$.  Therefore we can assume that
$$f_k(B) \cap Y'' = Z = C \times S^1 \times K \times I.
\eqno{(12)}
$$
Note that the volume of $B$ for our metrics is at most linear in $j$
by property (d) of Step 2 above.  Denote by $g_j(2)$ the resulting
modification of the metric $g_j(1)$.

\medskip\noindent
{\bf Step 4: Uniform lower bound for the $(k-1)$-systole.}  Let
$\alpha$ be a $(k-1)$-cycle in $X$ representing a class of infinite
order.  We replace $\alpha$ by a minimizing rectifiable current.

By Poincare duality, the cycle $\alpha$ must meet one of the manifolds
$B$ of step 3 (it is here that we need the assumption of Theorem 1.1
that $H_{k-1}(X)$ is torsion free).  By construction, the metrics
$g_j(2)$ of $X$ have bounded geometry (in the sense of 4.2) in a
neighborhood of $B$ of a fixed size.  Here the geometry near $B$ is
bounded inside $U_k$ by step~3, while the geometry near
$B\cap(X\setminus U_k)$ is bounded by step 1.  A lower bound for the
volume of $\alpha$ now follows from the monotonicity formula applied
at a point of the intersection $(\alpha\cap B)\subset X$.

\medskip\noindent
{\bf Step 5: The $(k-1,n-k+1)$-freedom.}  Now choose a basis
$C_{k-1,i}$ for $H_{k-1}(X)$, as in Lemma 3.1.  Let $U_{k-1}\subset X$
be a neighborhood, disjoint from $U_k$, of the $C_{k-1}$.  We modify
the metrics $g_j(2)$ in the neighborhood $U_{k-1}$, to construct
$(k-1,n-k+1)$-free metrics $g_j(3)$ with $\sys_{n-k+1} (g_j(3))=j^2$,
similarly to step 2.  Note that the volume of the manifold $B_{n-k+1}$
is increased to quadratic growth in $j$.

The metrics $g_j(3)$ cannot be chosen to dominate the metrics
$g_j(2)$.  Therefore we need to explain why our construction of
$(k-1,n-k+1)$-free metrics in the neighborhood $U_{k-1}\subset X$
respects the lower bound of step 4 for the $(k-1)$-systole.

We may assume that the manifold $B_{n-k+1}$ meets $C_{k-1}$
transversely in a finite number of points $x_1,\ldots, x_N$.  We
choose a neighborhood ${\cal O}\subset C_{k-1,i}$, such that $x_i\in
{\cal O}$ and ${\cal O}$ is diffeomorphic to a $(k-1)$-ball.

We choose a degree 1 map $C_{k-1,i}\to S^{k-1}$ to be a diffeomorphism
from ${\cal O}$ onto the complement of a point of the sphere, and
define the map $f_{k-1}: U_{k-1}\to S^{k-1}\times D^{n-k+1}$ as in
formula $(7)$.  Then we can ensure that, for each connected component
$(B_{n-k+1}\cap U_{k-1})_0$ of the intersection $B_{n-k+1}\cap
U_{k-1}$, we have
$$
f_{k-1}((B_{n-k+1}\cap U_{k-1})_0)= D^{n-k+1},
$$
so that again the metric near $B$ is bounded in all three regions,
$U_k$, $U_{n-k+1}$, and their complement in $X$.  The lower bound for
the $(k-1)$-systole now results from the monotonicity formula as
before.  This completes the proof of Theorem 6.1.

\medskip\noindent
The comments below were kindly provided by Shmuel Weinberger, for the
benefit of the reader who is a Riemannian geometer and is not
necessarily familiar with the details of the techniques of [T].

\medskip\noindent
{\bf Theorem 6.2} [T] Let $X$ be a compact manifold, $A$ a homology
class.  Then a multiple of $A$ can be represented by an imbedded
submanifold.

\medskip\noindent
{\it Idea of proof}.  If the codimension of $A$ in $X$ is odd, then a
multiple of $A$ can always be represented by a submanifold with
trivial normal bundle.  Indeed, the Poincare dual cohomology class is
a map to $K(\Z, odd)$, which is, at least rationally, $S^{odd}$.  Thus
a multiple is represented by a map to the sphere.  To construct the
desired submanifold, simply take the transverse inverse image of a
point of the sphere $S^{odd}$.

For even codimension $2n$ the proof is slightly more involved. 
We consider first an example in codimension 2.
Let $z\in H^2(X,\Z)$ be the dual 2-class to the homology class of
codimension 2 that we want to represent by a submanifold.  Let $f_z: X
\to \C P^{N}$ be the associated map to a skeleton of the classifying
space, for a sufficiently high $N$.  Then the transverse inverse image
of a hyperplane,
$$ f_z^{-1} (f_z(X) \cap \C P^{N-1}),
$$
is the desired submanifold.

Let $\zeta\to K(\Z,2)$ be the universal 2-bundle.  Its Euler class is
the generator of cohomology.  The Thom space $MG$ of $\zeta$ is again
the Eilenberg Maclane space: $MG = K(\Z,2)$.  Namely, $CP^N$ is the
one point compactification of $\zeta$ over $CP^{N-1}$.

Returning to the general case of even codimension $2n$, recall that
the cohomology of the Eilenberg Maclane space
$$H^*(K(Z,2n))\otimes \Q
$$ 
is a polynomial algebra on a $2n$-dimensional generator.  Now $BU(n)$
admits a map to a product  
$$
BU(n)\to K(\Z,2) \times K(\Z,4)\times \cdots \times K(\Z,2n),
$$  
which is defined (integrally) via Chern classes, and induces
isomorphism in rational cohomology.  

Now if $X$ and $Y$ are rationally equivalent, there may not in general
exist maps from $X$ to $Y$ and from $Y$ to $X$.  But for any fixed
$k$, one can map the $k$-skeleton of $X$ to $Y$ in such a way that the
map will induce an isomorphism in rational homology up to dimension
$n-1$.  This applies, for example, if $X$ is a finite cell complex,
and allows one to replace rational complexes by finite integral ones
in the argument that follows.

Thus there is an $n$-dimensional complex bundle over $K(\Z,2n)$ whose
$n$-th Chern class (\ie its Euler class) is a nonzero multiple of the
generator of
$$H^{2n}(K(\Z,2n)).
$$  
A Gysin sequence argument then shows that the Thom space $T$ of this
bundle is rationally $K(\Z,2n)$.  Thus one can (after multiplying)
lift the map $X \to K(\Z,2n)$ to $T$.  The desired submanifold of $X$
is then the transverse inverse image of the 0-section of $\zeta\subset
T$.

Note that this argument really only works on skeleta, while $K(\Z,2n)$
doesn't have the bundle, but rather, it rationally has a bundle.  It
integrally gets a bundle on its skeleta.

%\vfill\eject
\bigskip\noindent
{\large 7.  Curves and surfaces in $n$-manifolds}

\medskip\noindent
Our goal is to prove $(k,n-k)$ systolic freedom simultaneously for
$k=1$ and $k=2$, for a manifold with $H_1(X)=\Z$ (\cf Theorem 7.4
below).  However, we start with the simpler case of abelian
fundamental group.

\medskip\noindent
{\bf Proposition 7.2.}  Let $X$ be an orientable $n$-manifold with
free abelian fundamental group.  Then $X$ admits metrics which are
simultaneously $(1,n-1)$-free and $(2,n-2)$-free.

\medskip\noindent
{\it Proof.}  We will prove a stronger statement, namely, that $X$
admits metrics $g_j$ satisfying the inequality
$$
vol(g_j)\left( {1 \over \injrad(g_j) \sys_{n-1}(g_j)}+ {1 \over
\injrad(g_j)^2 \sys_{n-2}(g_j)}\right) \leq {\Const \over j}
\eqno{(13)}
$$
where $\injrad$ is the injectivity radius (\cf equations $(1)$ and
$(3)$ above).

Since the fundamental group is abelian, every incompressible
orientable surface in $X$ is either a sphere or a torus.  We choose a
rational basis $C_2=\cup_i C_{2,i}$ for $H_2(X)$, where each surface
$C_{2,i}$ is a torus $T^2$ with trivial normal bundle (spherical
classes may be represented by tori by adding a handle).  Then a
neighborhood $U_2$ of $C_2$ has the form $U_2=C_2\times D^{n-2}$.

We write $T^2= S^{2-1}\times C$.  Here the circle $S^{2-1}$ will play
a role similar to that of the sphere $S^{k-1}$ in the simply connected
case (Theorem 6.1).  The circle $C$ is similar to the circle $C$ in
the simply connected case, namely the second factor of the submanifold
$S^{k-1} \times C$ representing the class $[S^k]\in H_k(Y)$ (\cf
Proposition 3.3).  The difference from the simply connected case is
that now $C$ may not be contractible in $X$.  Here the analogue of the
map $f_k$ of formula $(7)$ is the map
$$
f_2: U_2\to S^{2-1}\times C\times D^{n-k}
$$ 
which is a diffeomorphism on each connected component of $U_2\subset
X$.

We choose a fixed smooth background metric $g$ on $X$.  We modify $g$
in $U_2$ to build $(2,n-2)$-free metrics $g_j$ on $X$ as in the proof
of Theorem 6.1, with quadratic growth of the $(n-2)$-systole:
$sys_{n-2}(g_j)
\sim j^2$, and linear volume growth: $\vol_n(g_j) \sim j$.  The
metrics $g_j$ have injectivity radius and 1-systole which are
uniformly bounded from below, since $f_2$ is a diffeomorphism.

We continue by modifying the metrics $g_j$ in a neighborhood
$U_1\subset X$ of a system of loops $C_1= \cup_i C_{1,i}$ generating
$\pi_1(X)=H_1(X)$.  Here each connected component of $U_1$ is
diffeomorphic to $C\times D^{n-1}$.  We choose a neighborhood
$Y''\subset U_1$ of the form $Y''=C\times S^1\times K\times I$ and
construct $(1,n-1)$-free metrics as in section 4.  The resulting
metrics satisfy quadratic growth in the $(n-1)$-systole, and obey a
uniform lower bound for the injectivity radius.  This concludes the
proof of Theorem 7.2.

\medskip\noindent
{\bf Lemma 7.3.}  Let $R_g$ be a closed orientable surface.  Let
$C_\alpha\subset R_g$ be a finite family of disjoint imbedded loops,
where $\alpha\in$ some index set.  Then there exists a smooth map $f:
R_g\to S^2$ to the sphere $S^2$, of nonzero degree, which maps each
loop $C_\alpha$ onto the equator of $S^2$, while a neighborhood of
$C_\alpha$ maps diffeomorphically onto a neighborhood of the equator.

\medskip\noindent
{\it Proof.}  Denote by $\NP\in S^2$ and $\SP\in S^2$ respectively the
north and south poles of the sphere.  It is obvious that there exists
such a map to the sphere with two points identified,
$S^2/(\SP\sim\NP)$.

To construct the map to the sphere itself, note that we may assume
that the family of loops is ``maximal'' in the sense that the
complement of $\cup_\alpha C_\alpha$ in $R_g$ is a collection of open
surfaces of one of two types: either a ``pair of pants'' or a cylinder
(some of the loops may be isotopic).  We choose disjoint tubular
neighborhoods, called annuli, $A_\alpha\subset R_g$ of the $C_\alpha$.
We identify the boundary of the annulus $\partial A_\alpha =
C_\alpha\times S^0$ with the product of the loop with the pair of
points $S^0=\{s,n\}$.  Let
$$
f: A_\alpha \to S^2\setminus \{\SP,\NP\}
$$ 
be an orientation-preserving diffeomorphism onto the complement of the
poles, while $f(C_\alpha)$ is the equator, so that $f$ extends by
continuity to the two boundary loops as follows:
$$f(C_\alpha\times \{s\})=\SP\in S^2, \,\, f(C_\alpha\times
\{n\})=\NP\in S^2.
$$ 
The complement $R_g\setminus (\cup_\alpha A_\alpha)$ of the annuli is
a disjoint union of connected ``complementary regions'' each of which
is either a pair of pants or a cylinder.  Each boundary loop of a
complementary region has a marking $\epsilon=\SP$ or $\epsilon=\NP$,
depending on whether it is the boundary loop of the adjacent annulus
$A_\alpha$ which is sent to $\SP\in S^2$ or to $\NP\in S^2$.  If the
marking $\epsilon$ is the same for all boundary loops of such a
complementary region, we collapse the entire region to the
corresponding pole $\epsilon\in S^2$.  If the boundary loops do not
all have the same markings, there are three possibilities:

1. The complementary region is a cylinder whose boundary loops are
marked $\SP$ and $\NP$.  Then we map it to $S^2$ by an
orientation-preserving diffeomorphism in the interior of the cylinder,
while the boundary loops are sent respectively to $\SP\in S^2$ and
$\NP\in S^2$;

2. The complementary region is a pair of pants and its triple of
   boundary loops are marked $\SP,\SP$, and $\NP$.  Then we subdivide
   the pair of pants into two regions, by choosing a new circle
   parallel to the boundary loop which is marked $\NP$, and cutting
   along it.  We label this new circle $\SP$ (see figure below).  Now,
   one of the regions is a cylinder, which can be handled in case 1.
   The other region is a pair of pants with the identical marking,
   $\SP$, on all three boundary loops, and we collapse it to the pole
   $\SP\in S^2$;

\medskip
\medskip\noindent
\centerline{{   		%%%the figure
\epsfxsize=6cm
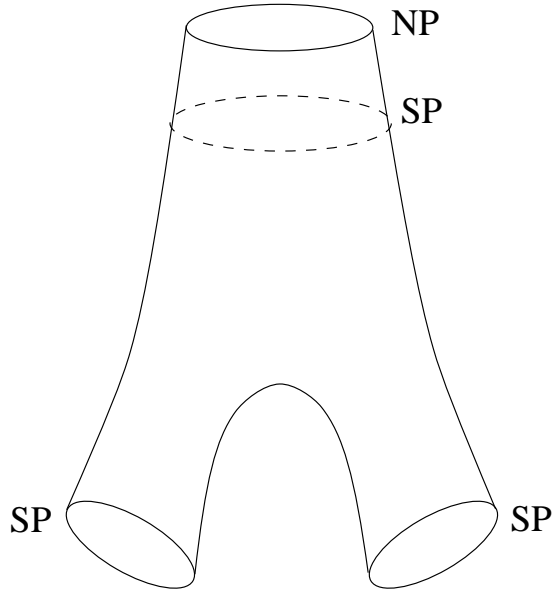
\epsfbox{geom.ps}  
}}

\medskip\noindent
3. The complementary region is a pair of pants and its boundary loops
   are marked $\NP,\NP$, and $\SP$.  This case is similar to the
   previous one.

\medskip\noindent
Since $f$ preserves orientation on each cylinder and pair of pants,
its degree is greater than or equal to the total number of circles
$C_\alpha$ at the outset, and is therefore nonzero (\cf calculation
following formula $(7)$).

\medskip\noindent
{\bf Theorem 7.4.}  Let $X$ be an orientable $n$-manifold.  Assume
that $H_1(X)=\Z$.  Then $X$ admits metrics which are simultaneously
$(1,n-1)$-free and $(2,n-2)$-free.

\medskip\noindent
{\it Proof.}  Let $C_2\subset X$ be a surface whose connected
components $C_{2,i}$ define a rational basis for 2-dimensional
homology.  Let $B=B_{n-1}$ be a hypersurface representing a generator
of the group $H_{n-1}(X)=\Z$.  We place $C_2$ and $B$ in transverse
position.  The intersection $C_2\cap B$ is a disjoint union of
imbedded circles $C_\alpha\subset C_2$.

The reason for requiring a unit first Betti number is that, for
$b_1\geq 2$, distinct components of $B_{n-1}$ may have a nonempty
triple intersection
$$
B_{n-1,i}\cap B_{n-1,j}\cap C_2\not= \emptyset,
\eqno{(14)}
$$ 
which by transversality is a finite set.  Thus, the collection of
loops on $C_2$ is in general not disjoint when $b_2\geq 2$, and the
construction of $(2,n-2)$-freedom cannot be accomplished with bounded
geometry near $B$.

Let $f:C_2\to S^2$ be the map of Lemma 7.3.  We now proceed as in the
proof of Theorem~6.1.  We first construct $(2,n-2)$-free metrics, as
in Proposition 5.2, by a modification in a neighborhood $U_2\subset X$
of the surface $C_2$, with bounded geometry in a neighborhood of the
loops $C_\alpha$.  Here the map of formula $(7)$ is replaced by the
map $f_2: C_2 \times D^{n-2}\to S^2\times D^{n-2}$, which is the map
$f$ of Lemma 7.3 times the identity on the disk.  Let $Y''\subset
U_{2,i}$ denote the submanifold $Y''= S^{2-1}\times C\times S^1\times
K\times I$ constructed as in section 3.  Then the hypersurface $B$ may
be assumed to satisfy
$$
f_2(B)\cap Y''= \cup_\alpha C_\alpha \times S^1 \times K \times I,
$$
so as to guarantee bounded geometry in a neighborhood of the
hypersurface $B$.

The resulting metrics may have zero injectivity radius, since the map
$f_2$ is no longer regular (as it was in Proposition 7.2).  However, a
loop representing a homology class of infinite order must meet the
hypersurface $B$ by Poincare duality.  Hence its length is bounded
below by the injectivity radius in a neighborhood of $B$.  Finally, we
modify the metric in a neighborhood of a loop $C_1$ representing a
generator of $H_1(X)$, as in Proposition 5.2, to attain
$(1,n-1)$-freedom.

\bigskip\noindent
{\large Acknowledgment.}  The author is grateful to I. Babenko for
catching an error in an earlier version of the paper, for signaling
the reference [W] (\cf 5.1 above), and for a number of helpful
suggestions.  Sh.~Weinberger contributed the comments following
Theorem 6.2, for which I thank him warmly.  The author expresses
appreciation to L.~Ambrosio, V.~Bangert, F.~Morgan, and B. White for
insightful comments, and to D. Garber for an expert drawing.
\vfill\eject
%\vskip1in\bigskip\noindent
\bigskip\noindent
{\large Bibliography}

\medskip\noindent
[B1] I.~Babenko, {\it Asymptotic invariants of smooth manifolds},
Russian Acad.\ Sci.\ Izv.\ Math.\ {41} (1993), 1--38.
\medskip\noindent
[B2] I.~Babenko, {\it Strong intersystolic softness of closed
manifolds}, Russian Math.\ Surveys, {55}, 5 (2000), 987--988.

\medskip\noindent
[B3] I. Babenko, Forte souplesse intersystolique de vari\'et\'es
ferm\'ees et de poly\`edres (2000), Universit\'e Montpellier-2,
pr\'epublication 13 (2000).

\medskip\noindent
[BK] I.~Babenko and M.~Katz, {\it Systolic freedom of
orientable manifolds}, Ann. Sci. \'{E}cole Norm. Sup. {31}
(1998), 787--809.

\medskip\noindent
[BKS] I.~Babenko, M.~Katz, and A.~Suciu, {\it Volumes,
middle-dimensional systoles, and Whitehead products},
Math. Res. Lett. {5} (1998), 461--471.  

\medskip\noindent
[BanK] V. Bangert and M. Katz, Stable systolic inequalities and
cohomology products, Communications on Pure and Applied Mathematics,
2002 (to appear).  

\medskip\noindent
[BeK] L.~Berard Bergery and M.~Katz, {\it Intersystolic inequalities
in dimension $3$}, Geom.\ Funct.\ Anal.\ {4} (1994), 621--632.

\medskip\noindent
[Be1] M. Berger, A l'ombre de Loewner, Ann.\ Scient.\ Ec.\ Norm.\ Sup.\
(1972) 241-260.

\medskip\noindent
[Be2] M. Berger, {\it Systoles et applications selon Gromov},
S\'eminaire N.~Bourbaki, expos\'{e} 771, Ast\'{e}risque {216} (1993),
279--310.

\medskip\noindent
[Fe1] H. Federer, Geometric measure theory.  Springer, 1969.

\medskip\noindent
[Fe2] H. Federer, Real flat chains, cochains, and variational problems,
Indiana Univ.\ Math.\ J.\ 24 (1974), 351-407.

\medskip\noindent
[Fr] M.~Freedman, {\it $\Z_2$ Systolic Freedom}, Geometry and Topology
Monographs, vol.~2 (1999), Proceedings of the Kirbyfest, paper no.~6,
pp.~113--123.

\noindent See {
http://www.maths.warwick.ac.uk/gt/GTMon2/paper6.abs.html}.

\medskip\noindent
[G1] M.~Gromov, {\it Filling Riemannian manifolds}, J. Differential
Geom.\ {18} (1983), 1--147.

\medskip\noindent
[G2] M. Gromov, Systoles and intersystolic inequalities, in Actes de
la Table Ronde de G\'{e}om\'{e}trie Diff\'{e}rentielle (Luminy, 1992),
S\'{e}min. Congr., vol.~1, Soc. Math. France, Paris, 1996,
pp.~291--362.

\medskip\noindent
[G3] M. Gromov, {\it Metric structures for Riemannian and
non-Riemannian spaces}.  Progress in Math., vol.~152, Birkh\"{a}user,
Boston, MA, 1999.
\medskip\noindent
[H] J.~Hebda, {\it The collars of a Riemannian manifold and stable
isosystolic inequalities}, Pacific J. Math.  {121} (1986), 339--356.
\medskip \noindent
[K1] M. Katz, Counterexamples to isosystolic inequalities, {\it
Geometriae Dedicata\/} 57 (1995) 195-206.
\medskip \noindent
[K2] M. Katz, Four-manifold systoles and surjectivity of period map,
preprint (2001).
\medskip\noindent
[KS1] M.~Katz and A.~Suciu, {\it Volume of Riemannian manifolds,
geometric inequalities, and homotopy theory}, in {\it Tel Aviv
Topology Conference: Rothenberg Festschrift} (M.~Farber, W.~L\"{u}ck,
and S.~Weinberger, eds.), Contemp. Math., vol.~231, Amer. Math. Soc.,
Providence, RI, 1999, pp.~113--136.
\medskip\noindent
[KS2] M. Katz and A. Suciu, Systolic freeedom of loop space, {\it
Geometric and Functional Analysis\/} 11 (2001) 60-73.

\medskip\noindent
[Mo] F. Morgan, {\it Geometric measure theory.  A beginner's guide}.
Academic Press, 1995.

\medskip\noindent
[P] C.~Pittet, {\it Systoles on $S^1\times S^n$}, Differential 
Geom.\ Appl.\ {7} (1997), 139--142.

\medskip\noindent
[T] R. Thom, Quelques proprietes globales des varietes
differentiables, Comment. Math. Helv. 28 (1954) 17-86.

\medskip\noindent
[W] Wright, Alden H.  Monotone mappings and degree one mappings
between $PL$ manifolds.  Geometric topology (Proc. Conf., Park City,
Utah, 1974), pp. 441--459.  Lecture Notes in Math., Vol. 438,
Springer, Berlin, 1975.

\vfill
\bigskip\noindent
Address: {Department of Mathematics and Statistics, Bar Ilan
University, Ramat Gan, 52900, Israel; and Departement de
Mathematiques, Universite Henri Poincare, Vandoeuvre 54506, France}.
Email address: {katzmik@macs.biu.ac.il}

\eject\end